\documentclass[journal]{IEEEtran}
\AtBeginDvi{} 

\IEEEoverridecommandlockouts

\usepackage[pdftex]{graphicx}

\usepackage{cite}
\usepackage{bm}
\usepackage{comment}
\usepackage{framed}
\usepackage{algorithm}
\usepackage{algorithmic}
\usepackage{slashbox}

\newtheorem{remark}{Remark}

\usepackage{ascmac}

\usepackage{amsmath}
\title{\LARGE \bf Riemannian optimal model reduction of stable linear systems}
\author{Kazuhiro~Sato
\thanks{K. Sato is with the School of Regional Innovation and Social Design Engineering,
 Kitami Institute of Technology,
 Hokkaido 090-8507, Japan,
email: ksato@mail.kitami-it.ac.jp}
}

\begin{document}

\maketitle

\begin{abstract}
In this paper, we develop a method for solving the problem of minimizing the $H^2$ error norm between the transfer functions of original and reduced systems on
the set of stable matrices and two Euclidean spaces.
That is, we develop  a method for identifying the optimal reduced system from all stable linear systems.
However, it is difficult to develop an algorithm for solving this problem, because the set of stable matrices is highly non-convex.
To overcome this issue, we show that the problem can be transformed into a tractable Riemannian optimization on 
 the product manifold of the set of skew-symmetric matrices, the manifold of the symmetric positive-definite matrices, and two Euclidean spaces.
The stability of the reduced systems constructed using the optimal solutions to our problem is preserved.
To solve the reduced problem, the Riemannian gradient and Hessian are derived and a Riemannian trust-region method is developed.
The initial point in the proposed approach is selected using the output from the balanced truncation (BT) method.
Numerical experiments demonstrate that our method considerably improves the results given by BT in the sense of the $H^2$ norm,
and also provides reduced systems that are globally near-optimal solutions to the problem of minimizing the $H^{\infty}$ error norm.  
Moreover, we show that our method provides a better reduced model than BT from the viewpoint of the frequency response.
\end{abstract}


%
\IEEEpeerreviewmaketitle

\section{Introduction}
%
%
%
%

Accurate modeling is essential to various system control methods.
However, the complexity of the controller is usually the same as that of the system.
That is, as the scale of the system to be controlled increases, the controller becomes more complex. 
This additional complexity can result in storage, accuracy, and computational speed problems \cite{antoulas2005approximation}.
Thus, we frequently need to approximate the original system as a small-scale model with high accuracy.

To produce a highly accurate reduced model, we use model reduction methods.
The most famous approach is the balanced truncation (BT) method \cite{antoulas2005approximation, dullerud2000course, moore1981principal}.
BT provides a stable reduced model with guaranteed $H^{\infty}$ bounds, as long as the original model is stable.
Another famous technique is the moment matching method \cite{antoulas2010interpolatory, astolfi2010model, gugercin2008h_2, ionescu2014families},
which produces a reduced system matching some coefficients of the transfer function of a given linear system.
In \cite{sato2015riemannian,  yan1999approximate},
the $H^2$ optimal model reduction problem was studied for general stable linear systems by formulating the optimization problem on the Stiefel manifold.
However, the methods developed in \cite{sato2015riemannian,  yan1999approximate} could be improved further, because
they only search for the optimal reduced model from a subset of all stable linear systems.

In this study, we develop a novel $H^2$ optimal model reduction method for stable linear systems.
The problem is formulated as a minimization problem of the $H^2$ error norm between the transfer functions of the original and reduced systems on the product set of stable matrices and two Euclidean spaces.
That is, unlike \cite{sato2015riemannian,  yan1999approximate}, we search for the optimal reduced model with respect to all stable linear systems. 
However, it is difficult to develop an algorithm for solving this problem, because the set of stable matrices is highly non-convex \cite{orbandexivry2013nearest}.

The contributions of this paper can be summarized as follows.
We show that the original difficult problem can be transformed into a tractable Riemannian optimization problem on the product manifold of the vector space of skew symmetric matrices, the manifold of symmetric positive-definite matrices, and two Euclidean spaces.
Thus, we propose a Riemannian trust-region method for solving the model reduction problem.
To this end, we derive the Riemannian gradient and Hessian of the objective function.
The initial point is given by the result of the BT method.
Numerical experiments demonstrate that our proposed method improves the results of the BT method in the sense of the $H^2$ and $H^{\infty}$ norms.
That is, although the aim of our optimization problem is to minimize the $H^2$ error norm between the transfer functions of the original and reduced systems,
the $H^{\infty}$ error norm between those is also smaller than that of the BT method.
Furthermore, we illustrate that our proposed method produces reduced systems that are globally near-optimal solutions to the problem of minimizing the $H^{\infty}$ error norm.
Moreover, we show that our method provides a better reduced model than the BT method from the viewpoint of the frequency response.

The remainder of this paper is organized as follows.
In Section \ref{sec2}, we formulate the $H^2$ optimal model reduction problem on the set of stable matrices and two Euclidean spaces.
In Section \ref{sec3}, we transform the problem into a tractable Riemannian optimization problem. 
In Section \ref{sec4}, we propose an optimization algorithm for solving our problem and a technique for choosing the initial point.
In Section \ref{sec5}, we demonstrate that our method is more effective than the BT method when the dimension of the reduced system is small.
Finally, our conclusions are presented in Section \ref{sec6}.

{\it Notation:} The sets of real and complex numbers are denoted by ${\bf R}$ and ${\bf C}$, respectively.
The identity matrix of size $n$ is denoted by $I_n$.
The symbol ${\rm Skew}(n)$ denotes the vector space of skew-symmetric matrices in ${\bf R}^{n\times n}$.
The manifold of  symmetric positive-definite matrices in ${\bf R}^{n\times n}$ is denoted by ${\rm Sym}_+(n)$.
The tangent space at $x$ on a manifold $X$ is denoted by $T_x X$.
Given a matrix $A\in {\bf R}^{n\times n}$, ${\rm tr} (A)$ denotes the sum of the elements on the diagonal of $A$, and $A_{i,j}$ denotes the entry in row $i$ and column $j$.
Moreover, ${\rm sym}(A)$ and ${\rm sk}(A)$ denote the symmetric and skew-symmetric parts of $A$, respectively; i.e., ${\rm sym}(A)=\frac{A+A^T}{2}$ and ${\rm sk}(A)=\frac{A-A^T}{2}$.
Here, $A^T$ denotes the transpose of $A$.
Given a vector $v\in {\bf C}^n$, $||v||$ denotes the Euclidean norm.
The Hilbert space $L^2({\bf R}^n)$ is defined by
$L^2({\bf R}^n) := \left\{ f: [0,\infty) \rightarrow {\bf R}^n\, \big|\, \int_0^{\infty} ||f(t)||^2 dt <\infty \right\}$.
Given a measurable function $f: [0,\infty)\rightarrow {\bf R}^n$, $||f||_{L^2}$ and $||f||_{L^{\infty}}$ denote the $L^2$ and $L^{\infty}$ norms of $f$, respectively, i.e.,
\begin{align*}
& ||f||_{L^2}:= \sqrt{ \int_0^{\infty} ||f(t)||^2 dt }, \\
& ||f||_{L^{\infty}} := \sup_{t\geq 0} ||f(t)||. 
\end{align*}
Given a matrix $A\in {\bf C}^{n\times n}$, $||A||$ and $||A||_F$ denote the induced and Frobenius norms, i.e.,
\begin{align*}
||A|| &:= \sup_{v\in C^n\backslash \{0\}} \frac{ ||Av||}{||v||}, \\
||A||_F &:= \sqrt{ {\rm tr} (A^* A)},
\end{align*}
where the superscript $*$ denotes  Hermitian conjugation, and
${\rm tr} (A)$ is the trace of $A$, i.e., the sum of the diagonal elements of $A$.
For a matrix function $G(s)\in {\bf C}^{n\times n}$, $||G||_{H^2}$ and $||G||_{H^{\infty}}$ denote the $H^2$ and $H^{\infty}$ norms of $G$, respectively, i.e.,
\begin{align*}
& ||G||_{H^2}  := \sqrt{ \frac{1}{2\pi} \int_{-\infty}^{\infty} ||G(i\omega)||_F^2 d \omega}, \\
& ||G||_{H^{\infty}} := \sup_{\omega \in{\bf R}} \bar{\sigma} (G(i\omega)),
\end{align*}
where $i$ is the imaginary unit, and $\bar{\sigma}(G(i\omega))$ denotes the maximum singular value of $G(i\omega)$. 

\section{Problem Setup} \label{sec2}

This section describes the formulation of our problem.

As the original system, we consider the linear continuous-time system
\begin{align}
\begin{cases}
\dot{x} = Ax + Bu, \\
y = C x,
\end{cases} \label{1}
\end{align}
where $x\in {\bf R}^{n}$, $u\in {\bf R}^{m}$, and $y\in {\bf R}^{p}$ are the state, input, and output, respectively.
The matrices $A\in {\bf R}^{n\times n}$, $B\in {\bf R}^{n\times m}$, and $C\in {\bf R}^{p\times n}$ are constant matrices.
Throughout this paper, we assume that system \eqref{1} is asymptotically stable; i.e., the real parts of all the eigenvalues of the matrix $A$ are negative.
We also call the matrix $A$ stable if system \eqref{1} is asymptotically stable.

In this paper, we consider the following $H^2$ optimal model reduction problem of preserving the stability.
\begin{framed}
Problem 1:
\begin{align*}
&{\rm minimize} \quad ||G-\hat{G}_r||_{H^2} \\
&{\rm subject\, to} \quad (A_r,B_r,C_r)\in {\bf S}^{r\times r} \times {\bf R}^{r\times m}\times {\bf R}^{p\times r}.
\end{align*}
\end{framed}

\noindent
Here, $G$ is the transfer function of system \eqref{1}, i.e.,
\begin{align*}
G(s) := C(sI_n -A)^{-1} B,\quad s\in {\bf C},
\end{align*}
$\hat{G}_r$ is the transfer function of the reduced system
\begin{align}
\begin{cases}
\dot{\hat{x}}_r = A_r\hat{x}_r + B_ru, \\
\hat{y}_r = C_r \hat{x}_r,
\end{cases} \label{remark_red}
\end{align}
and ${\bf S}^{r\times r}$ denotes the set of all stable matrices.
Note that if $u\in L^2({\bf R}^m)$, then the error $y-\hat{y}_r$ satisfies
\begin{align}
||y-\hat{y}_r||_{L^{\infty}} \leq ||G-\hat{G}_r||_{H^2} \cdot ||u||_{L^2}. \label{3_2}
\end{align}
The proof is shown in Appendix \ref{ape0}.
That is, if $||G-\hat{G}_r||_{H^2}$ is sufficiently small, then we can expect $||y-\hat{y}_r||_{L^{\infty}}$ to become almost zero for any $u$ with a small $||u||_{L^2}$.

It is difficult to solve Problem 1 because 
the set ${\bf S}^{r\times r}$ is highly non-convex \cite{orbandexivry2013nearest}.
To develop an algorithm for solving Problem 1, we transform Problem 1 into an equivalent, tractable Riemannian optimization problem.

\section{Equivalent Riemannian optimization problem} \label{sec3}

This section proves that Problem 1 is equivalent to 

\begin{framed}
Problem 2:
\begin{align*}
&{\rm minimize} \quad f(J_r,R_r,B_r,C_r):=||G-G_r||_{H^2}^2 \\
&{\rm subject\, to} \quad (J_r,R_r,B_r,C_r)\in M.
\end{align*}
\end{framed}

\noindent
Here, 
$G_r$ is the transfer function of the reduced system
\begin{align}
\begin{cases}
\dot{x}_r = (J_r-R_r)x_r + B_r u, \\
y_r = C_r x_r,
\end{cases} \label{reduction} 
\end{align}
and
\begin{align*}
M := {\rm Skew}(r)\times {\rm Sym}_+(r) \times {\bf R}^{r\times m} \times {\bf R}^{p\times r}.
\end{align*}
Note that system \eqref{reduction} is also asymptotically stable, because the real parts of all the eigenvalues of $J_r-R_r$ are negative.

To this end, we first note that,
because it is asymptotically stable, system \eqref{1} can be transformed into 
\begin{align}
\begin{cases}
\dot{x} = (J-R)\mathcal{Q}x + Bu, \\
y = C x,
\end{cases} \label{11}
\end{align}
where $\mathcal{Q}\in {\rm Sym_+}(n)$ and
\begin{align*}
J &:= \frac{1}{2} (A\mathcal{Q}^{-1}-\mathcal{Q}^{-1}A^T)\in {\rm Skew}(n), \\
R &:=  -\frac{1}{2} (A\mathcal{Q}^{-1}+\mathcal{Q}^{-1}A^T)\in {\rm Sym}_+(n).
\end{align*}
Although the proof can be found in Proposition 1 of \cite{prajna2002lmi}, we repeat it here for completeness. 
From the asymptotic stability of system \eqref{1}, it follows that
there exists a Lyapunov function of the form $\mathcal{V}(x) = \frac{1}{2} x^T \mathcal{Q} x$ with $\mathcal{Q}\in {\rm Sym_+}(n)$  such that $\dot{\mathcal{V}}(x)= \frac{1}{2} x^T(A^T\mathcal{Q}+\mathcal{Q}A)x<0$, i.e., 
\begin{align}
-(A^T\mathcal{Q}+\mathcal{Q}A)\in {\rm Sym}_+(n). \label{hogehogehoge}
\end{align}
Thus, $J\in {\rm Skew}(n)$, $R\in {\rm Sym}_+(n)$, and $(J-R)\mathcal{Q}=A$.
Note that we can easily find $\mathcal{Q}\in {\rm Sym}_+(n)$ satisfying \eqref{hogehogehoge}.
In fact, because the matrix $A$ is stable, there exists $\mathcal{Q}\in {\rm Sym}_+(n)$ satisfying the Lyapunov equation
\begin{align}
A^T\mathcal{Q}+\mathcal{Q}A + I_n = 0, \label{lyap1}
\end{align}
as shown in \cite{dullerud2000course}.
The Lyapunov equation in \eqref{lyap1} can be efficiently solved using the Bartels--Stewart algorithm \cite{bartels1972solution}.
Because the transfer function of \eqref{11} coincides with that of \eqref{1},
 Problem 1 is equivalent to 

\begin{framed}
Problem 3:
\begin{align*}
&{\rm minimize} \quad ||G-\check{G}_r||_{H^2} \\
&{\rm subject\, to} \quad (J_r,R_r,\mathcal{Q}_r,B_r,C_r)\in N.
\end{align*}
\end{framed}

\noindent
Here, $\check{G}_r$ is the transfer function of the reduced system
\begin{align*}
\begin{cases}
\dot{\check{x}}_r = (J_r-R_r)\mathcal{Q}_r\check{x}_r + B_ru, \\
\check{y}_r = C_r \check{x}_r,
\end{cases} 
\end{align*}
and $N:= {\rm Skew}(r)\times {\rm Sym}_+(r) \times {\rm Sym}_+(r)\times {\bf R}^{r\times m}\times {\bf R}^{p\times r}$.

Next, we show that Problem 3 can be transformed into Problem 2.
To see this, we note that system \eqref{11} is equivalent to the form
\begin{align}
\begin{cases}
\dot{\tilde{x}} = (\tilde{J}-\tilde{R})\tilde{x} + \tilde{B}u, \\
y = \tilde{C} \tilde{x},
\end{cases} \label{original}
\end{align}
where $\tilde{x}\in{\bf R}^n$, $\tilde{J}\in {\rm Skew}(n)$, $\tilde{R}\in {\rm Sym}_+(n)$, $\tilde{B}\in{\bf R}^{n\times m}$, and $\tilde{C}\in {\bf R}^{p\times n}$.
In fact, because $\mathcal{Q}$ is a positive-symmetric matrix, there exists a unique lower triangular $L \in {\bf R}^{n\times n}$ with positive diagonal entries such that $\mathcal{Q} = LL^T$. This is called the Cholesky decomposition of $\mathcal{Q}$. For a detailed explanation, see \cite{golub2012matrix}.
Thus, if we perform a coordinate transformation $\tilde{x} = (L^{-1})^Tx$,
we obtain \eqref{original}, where $\tilde{J} = L^T JL$, $\tilde{R}=L^T RL$, $\tilde{B}=L^TB$, and $\tilde{C}=C(L^{-1})^T$. 
Because the transfer function of \eqref{original} coincides with that of \eqref{11}, Problem 3 is equivalent to Problem 2.

From the above discussion, Problem 2 is equivalent to Problem 1, which completes the proof.

In contrast to Problem 1, we can develop an algorithm for solving Problem 2 using a Riemannian optimization method \cite{absil2009optimization}, as shown in the next section.

\begin{remark}
Reference \cite{van2008h2} considered
\begin{align*}
&{\rm minimize} \quad ||G-\hat{G}_r||_{H^2}^2 \\
&{\rm subject\, to} \quad (\hat{A}_r,\hat{B}_r,\hat{C}_r)\in {\bf R}^{r\times r} \times {\bf R}^{r\times m}\times {\bf R}^{p\times r},
\end{align*}
and proved that if reduced system \eqref{remark_red} is controllable and observable, then at every stationary point of $||G-\hat{G}_r||_{H^2}^2$, we have that
\begin{align*}
& \hat{A}_r = W^TAV,\,\, \hat{B}_r = W^TB,\,\, \hat{C}_r = CV,\,\, W^TV=I_r.
\end{align*}
Based on this fact, \cite{gugercin2008h_2, antoulas2010interpolatory} developed an algorithm for finding such $V$ and $W$.
Although the algorithm can be applied to the model reduction of large-scale systems, a sequence produced by the algorithm does not generally converge to a local optimal solution, except for single-input--single-output symmetric systems \cite{flagg2012convergence}.
\end{remark}

\begin{remark}
References \cite{sato2015riemannian, yan1999approximate} considered
\begin{framed}
Problem 4:
\begin{align*}
&{\rm minimize} \quad ||G-\bar{G}_r||_{H^2}^2 \\
&{\rm subject\, to} \quad U\in {\rm St}(r,n).
\end{align*}
\end{framed}

\noindent
Here, $\bar{G}_r$ is the transfer function of the reduced system
\begin{align*}
\begin{cases}
\dot{\hat{x}}_r = U^TAU\bar{x}_r + U^TBu, \\
\bar{y}_r = CU \bar{x}_r,
\end{cases}
\end{align*}
and ${\rm St}(r,n)$ is the Stiefel manifold defined by
\begin{align*}
{\rm St}(r,n):= \{ U\in {\bf R}^{n\times r}\,|\, U^TU = I_r\}.
\end{align*}
As explained in \cite{sato2015riemannian, yan1999approximate}, if $A+A^T$ is negative-definite,
then $A$ and $U^TAU$ are stable, i.e., $A\in {\bf S}^{n\times n}$ and $U^TAU\in {\bf S}^{r\times r}$.
Thus, if this is the case, a solution to Problem 4 is a feasible solution to Problem 1.
That is, by solving Problem 4, we can obtain feasible solutions to Problem 1.
However, in general, the optimal value of Problem 4 is larger than that of Problem 1 \cite{sato2018automatica, sato2018TAC}.
This is because any method for solving Problem 4 searches for the optimal reduced system from a subset of the stable linear systems.
\end{remark}

\begin{remark}
Instead of Problem 2, we can consider the following $H^{\infty}$ optimal model reduction problem.
\begin{framed}
Problem 5:
\begin{align*}
&{\rm minimize} \quad ||G-G_r||_{H^{\infty}} \\
&{\rm subject\, to} \quad (J_r,R_r,B_r,C_r)\in M.
\end{align*}
\end{framed}

\noindent
However, in contrast to Problem 2, the objective function $||G-G_r||_{H^{\infty}}$ is not differentiable.
Thus, it is difficult to develop an algorithm for solving Problem 5.
In Section \ref{sec4}, we demonstrate that there are examples for which we can obtain a globally near-optimal solution to Problem 5 by solving Problem 2.
\end{remark}

\section{Optimization algorithm for Problem 2} \label{sec4}

\subsection{Riemannian gradient, Hessian, and exponential map} \label{sec4A}

To develop an optimization algorithm for solving Problem 2, we derive the Riemannian gradient and Hessian of the objective function $f$, and compute the exponential map on the manifold $M$.

To this end, we first note that,
because systems \eqref{original} and \eqref{reduction} are asymptotically stable, the objective function $f$ can be expressed as
\begin{align*}
f(J_r,R_r,B_r,C_r) &= {\rm tr}( \tilde{C}\Sigma_c \tilde{C}^T +C_r PC_r^T-2C_rX^T\tilde{C}^T)  \\
&= {\rm tr}( \tilde{B}^T\Sigma_o \tilde{B} +B_r^T Q B_r+2\tilde{B}^T Y B_r),
\end{align*}
where $\Sigma_c$, $\Sigma_o$, $P$, $Q$, $X$, and $Y$ satisfy
\begin{align}
(\tilde{J}-\tilde{R})\Sigma_c +\Sigma_c (\tilde{J}-\tilde{R})^T +\tilde{B}\tilde{B}^T  &= 0,  \nonumber\\
(\tilde{J}-\tilde{R})^T\Sigma_o +\Sigma_o (\tilde{J}-\tilde{R}) +\tilde{C}^T\tilde{C}  &= 0,  \nonumber\\
(J_r-R_r)P+P(J_r-R_r)^T+B_rB_r^T &=0, \label{3} \\
(J_r-R_r)^T Q+Q(J_r-R_r)+C_r^T C_r &=0, \label{4} \\
(\tilde{J}-\tilde{R})X+X(J_r-R_r)^T+\tilde{B}B_r^T &=0, \label{5} \\
(\tilde{J}-\tilde{R})^TY+Y(J_r-R_r)-\tilde{C}^TC_r &=0, \label{6}
\end{align} 
respectively.
For a detailed derivation, see \cite{sato2015riemannian, van2008h2, yan1999approximate}.

Let $\bar{f}$ denote the extension of the objective function $f$ to the Euclidean space ${\bf R}^{r\times r}\times {\bf R}^{r\times r} \times {\bf R}^{r\times m}\times {\bf R}^{p\times r}$.
In the same way as in previous studies \cite{sato2016new, sato2017, sato2017_2, sato2018TAC, sato2018automatica}, we then obtain 
\begin{align}
& \nabla \bar{f}(J_r,R_r,B_r,C_r) \nonumber \\
=& 2(QP+Y^TX, -QP-Y^TX, QB_r+Y^TB, C_rP-CX). \label{16}
\end{align}
To derive the Riemannian gradient and Hessian, we define the Riemannian metric of the manifold $M$ as
\begin{align}
& \langle (\xi_1,\eta_1,\zeta_1,\kappa_2),(\xi_2,\eta_2,\zeta_2,\kappa_2) \rangle_{(J_r,R_r,B_r,C_r)} \nonumber \\
:=&{\rm tr} (\xi_1^T\xi_2) + {\rm tr}( R_r^{-1} \eta_1 R_r^{-1} \eta_2  ) +  {\rm tr}(\zeta_1^T \zeta_2) +  {\rm tr}(\kappa_1^T \kappa_2) \label{Riemannian_metric}
\end{align}
for $(\xi_1,\eta_1,\zeta_1,\kappa_1),(\xi_2,\eta_2,\zeta_2,\kappa_2) \in T_{(J_r,R_r,B_r,C_r)} M$.
It then follows from \eqref{gradient} in Appendix \ref{apeB} and \eqref{16} that
\begin{align}
& {\rm grad}\, f(J_r,R_r,B_r,C_r) \nonumber\\
 =& (2{\rm sk}(QP+Y^TX), -2 R_r {\rm sym}(QP+Y^TX)R_r,  \nonumber \\
&\,\, 2(QB_r+Y^TB), 2(C_rP-CX)). \label{grad_f}
\end{align}
Furthermore, from \eqref{Hess} in Appendix \ref{apeB} and \eqref{16}, the Riemannian Hessian of $f$ at $(J_r,R_r,B_r,C_r)$ is given by
\begin{align}
&  {\rm Hess}\, f(J_r,R_r,B_r,C_r) [(J'_r,R'_r,B'_r,C'_r)] \nonumber \\
=& ( 2{\rm sk}(Q'P+QP'+Y'^TX+Y^TX'), \nonumber \\
& -2R_r {\rm sym}( Q'P+QP'+Y'^TX+Y^TX') R_r \nonumber \\
&\, -2 {\rm sym}(R'_r {\rm sym} (QP+Y^TX) R_r), \nonumber\\
&\,\, 2(Q'B_r+QB'_r+Y'^TB), 2(C'_rP+C_rP'-CX') ), \label{Hess_f}
\end{align}
where $P'$, $Q'$, $X'$, and $Y'$  are the solutions to
\begin{align*}
& (J_r-R_r)P' + P'(J_r-R_r)^T + (J'_r-R'_r)P + P(J'_r-R'_r)^T \\
&+ B'_rB_r^T + B_rB'^T_r = 0, \\
& (J_r-R_r)^T Q' +Q' (J_r-R_r)+(J'_r-R'_r)^T_r Q+Q (J'_r-R'_r) \\
&+ C'_r C_r^T+C_r C'^T_r =0,  \\
& (\tilde{J}-\tilde{R})^TX'+X'(J_r-R_r)+X(J'_r-R'_r)+\tilde{B}B'^T_r =0, \\
& (\tilde{J}-\tilde{R})^TY'+Y'(J_r-R_r)+Y(J'_r-R'_r)-\tilde{C}^TC'_r =0, 
\end{align*}
respectively. Note that these equations are obtained by differentiating \eqref{3}, \eqref{4}, \eqref{5}, and \eqref{6}, respectively.
Moreover, from \eqref{8} in Appendix A, we can define the exponential map on the manifold $M$ as
\begin{align}
& {\rm Exp}_{(J_r,R_r,B_r,C_r)}(\xi,\eta,\zeta,\kappa) \nonumber \\
:=& (J_r+\xi,{\rm Exp}_{R_r}(\eta), B_r+\zeta,C_r+\kappa) \nonumber \\
=& (J_r+\xi, R_r^{\frac{1}{2}} \exp (R_r^{-\frac{1}{2}} \eta R_r^{-\frac{1}{2}} ) R_r^{\frac{1}{2}}, B_r+\zeta,C_r+\kappa) \label{33}
\end{align}
for any $(\xi,\eta,\zeta,\kappa)\in T_{(J_r,R_r,B_r,C_r)} M$.

\subsection{Trust-region method for Problem 2}

Algorithm 1 describes the Riemannian trust-region method for solving Problem 2.
At each iterate $p_r:=(J_r,R_r,B_r,C_r)\in M$ in the Riemannian trust-region method, 
we evaluate the quadratic model $\hat{m}_{p_r}$ of the objective function $f$ within a trust region:
\begin{align*}
&\quad \hat{m}_{p_r}(\xi,\eta,\zeta,\kappa)  \\
=& f(J_r,R_r,B_r,C_r)  + \langle {\rm grad}\,f(J_r,R_r,B_r,C_r), (\xi,\eta,\zeta,\kappa) \rangle_{p_r} \nonumber\\
&+\frac{1}{2}  \langle {\rm Hess}\, f(J_r,R_r,B_r,C_r)[(\xi,\eta,\zeta,\kappa)], (\xi,\eta,\zeta,\kappa) \rangle_{p_r}.
\end{align*}
Because we can construct the gradient and Hessian of $f$ as in Section \ref{sec4A},
we can construct $\hat{m}_{p_r}$.
A trust region with a radius $\Delta>0$ at $p_r\in M$ is defined as a ball in $T_{p_r} M$.
The trust-region sub-problem at $p_r\in M$ with the radius $\Delta$ is thus defined as the problem of minimizing $\hat{m}_{p_r}(\xi,\eta,\zeta,\kappa)$ subject to
$(\xi,\eta,\zeta,\kappa)\in T_{p_r}M$, $||(\xi,\eta,\zeta,\kappa)||_{p_r}:= \sqrt{\langle (\xi,\eta,\zeta,\kappa),(\xi,\eta,\zeta,\kappa) \rangle_{p_r}}\leq \Delta$.
This sub-problem can be solved by the truncated conjugate gradient method \cite{absil2009optimization}.
We then compare the decrease in the objective function $f$ and the model $\hat{m}_{p_r}$ attained by the resulting $(\xi_*,\eta_*,\zeta_*,\kappa_*)$, and use this to determine
whether $(\xi_*,\eta_*,\zeta_*,\kappa_*)$ should be accepted and whether the trust region of radius $\Delta$ is appropriate.
The constants $1/4$ and $3/4$ in the conditional expressions in Algorithm 1 are commonly used in the trust-region method for a general unconstrained optimization problem.
These values ensure the convergence properties of the algorithm \cite{absil2009optimization}.
In fact, if the trust-region sub-problem is carefully solved, sequences generated by the Riemannian trust-region method converge quadratically under certain assumptions on the objective function in question \cite{absil2009optimization}.

Note that the reduced system attained by Algorithm \ref{algorithm} is asymptotically stable,
because $(J_r,R_r,B_r,C_r)\in M$ at each iteration.

\begin{algorithm}                      
\caption{Trust-region method for Problem 2.}    \label{algorithm}     
\label{alg1}                          
\begin{algorithmic}[1]
{\small
\STATE Choose an initial point $(p_r)_0 \in M$ and parameters $\bar{\Delta}>0$, $\Delta_0\in (0,\bar{\Delta})$, $\gamma'\in [0,\frac{1}{4})$.
\FOR{$k=0,1,2,\ldots$ }
\STATE Solve the following trust-region sub-problem for $(\xi,\eta,\zeta,\kappa)$ to obtain $(\xi_k,\eta_k,\zeta_k,\kappa_k)\in T_{(p_r)_k} M$:
\begin{align*}
&{\rm minimize}\quad \hat{m}_{(p_r)_k}(\xi,\eta,\zeta,\kappa) \\
&{\rm subject\, to}\quad ||(\xi,\eta,\zeta,\kappa)||_{(p_r)_k} \leq \Delta_k, \\
&{\rm where}\quad  (\xi,\eta,\zeta,\kappa)\in T_{(p_r)_k}M.
\end{align*}
\STATE Evaluate
\begin{align*}
&\gamma_k := \frac{ f({\rm Exp}_{(p_r)_k}(0,0,0,0)) -f({\rm Exp}_{(p_r)_k}(\xi_k,\eta_k,\zeta_k,\kappa_k))}{ \hat{m}_{(p_r)_k}(0,0,0,0)- \hat{m}_{(p_r)_k} (\xi_k,\eta_k,\zeta_k,\kappa_k)}.
\end{align*} 
\IF {$\gamma_k<\frac{1}{4}$}
\STATE 
$\Delta_{k+1}=\frac{1}{4}\Delta_k$.
\ELSIF {$\gamma_k>\frac{3}{4}$ and $||(\xi_k,\eta_k,\zeta_k,\kappa_k)||_{(p_r)_k} = \Delta_k$}
\STATE
$\Delta_{k+1} = \min (2\Delta_k,\bar{\Delta})$.
\ELSE 
\STATE
$\Delta_{k+1} = \Delta_k$.
\ENDIF
\IF {$\gamma_k>\gamma'$}
\STATE
$(p_r)_{k+1} = {\rm Exp}_{(p_r)_k}(\xi_k,\eta_k,\zeta_k,\kappa_k)$.
\ELSE
\STATE
$(p_r)_{k+1} = (p_r)_{k}$.
\ENDIF
\ENDFOR
}
\end{algorithmic}
\end{algorithm}

\subsection{Initial point in Algorithm 1}

In this subsection, we describe a technique for choosing the initial point $(p_r)_0\in M$ in Algorithm 1 using the output of the BT method \cite{antoulas2005approximation, dullerud2000course, moore1981principal}.
The BT method can be implemented using the MATLAB command {\it balred} (i.e., we can easily implement the BT method),
and provides satisfactory reduced models in many cases.
 
The BT method outputs the reduced matrices $(A_r)_{{\rm BT}}$, $(B_r)_{{\rm BT}}$, and $(C_r)_{{\rm BT}}$; the matrix $(A_r)_{{\rm BT}}$ is stable, because the original matrix $A$ is stable  \cite{antoulas2005approximation, dullerud2000course, moore1981principal}.
Thus, there exists $\mathcal{Q}_r\in {\rm Sym}_+(r)$ satisfying 
\begin{align*}
(A_r)_{{\rm BT}}^T\mathcal{Q}_r +\mathcal{Q}_r (A_r)_{{\rm BT}} + I_r = 0,
\end{align*}
as explained in Section \ref{sec3}.
Next, we define
\begin{align*}
(J_r)_{{\rm BT}} &:= \frac{1}{2} \left((A_r)_{{\rm BT}}\mathcal{Q}^{-1}_r-\mathcal{Q}^{-1}_r(A_r)^T_{{\rm BT}} \right), \\
(R_r)_{{\rm BT}} &:=  -\frac{1}{2} \left((A_r)_{{\rm BT}}\mathcal{Q}^{-1}_r+\mathcal{Q}^{-1}_r(A_r)^T_{{\rm BT}} \right).
\end{align*}
Finally, we perform the Cholesky decomposition of $\mathcal{Q}_r=L_rL_r^T$, and set the initial point
\begin{align*}
&(p_r)_0=((J_r)_0, (R_r)_0,(B_r)_0,(C_r)_0) \\
=& (L_r^T(J_r)_{{\rm BT}}L_r, L_r^T(R_r)_{{\rm BT}}L_r, L^T_r(B_r)_{{\rm BT}}, (C_r)_{{\rm BT}} (L_r^{-1})^T).
\end{align*}
Note that, because transfer functions are invariant under coordinate transformations, we have that
\begin{align*}
(G_r)_{{\rm BT}}&= (C_r)_{{\rm BT}}(sI_r -(A_r)_{{\rm BT}})^{-1}(B_r)_{{\rm BT}} \\
&= (C_r)_0(sI_r- ( (J_r)_0 - (R_r)_0))^{-1} (B_r)_0,
\end{align*}
where $(G_r)_{BT}$ is the transfer function of the reduced system attained by the BT method.

\section{Numerical Experiments} \label{sec5}

In this section, two examples are presented to illustrate that our method improves the BT result in the sense of the $H^2$ norm.
Furthermore, we show that our method may provide better results for the $H^{\infty}$ norm and the frequency response than the BT method.
To this end, we have used Manopt \cite{boumal2014manopt}, which is a MATLAB toolbox for optimization on manifolds.

\subsection{Mass-spring-damper system}

We consider mass-spring-damper systems with masses $m_i$, spring constants $k_i$, and damping constants $c_i$ $(i=1,2,\ldots,\frac{n}{2})$, where $n$ is an even number.
The inputs $u_1$ and $u_2$ are the external forces applied to the first two masses, $m_1$ and $m_2$.
The output $y_1$ is the displacement of mass $m_1$.
The state variables $\tilde{x}_j$ $(j=1,3,\ldots)$ are the displacements of mass $m_j$ and 
the state variables $\tilde{x}_k$ $(k=2,4,\ldots)$ are the momentums of mass $m_k$.
Here, we only consider the case where $m_i=4$, $k_i=4$, and $c_i=1$ $(i=1,2,\ldots,\frac{n}{2})$.
The system can be described by \eqref{11} and the system matrices are given by
$\tilde{J}_{1,2} =\tilde{J}_{3,4} = \cdots = \tilde{J}_{(n-1),n} = 1$,
$\tilde{J}_{2,1} = \tilde{J}_{4,3} = \cdots = \tilde{J}_{n,(n-1)} = -1$,
$\tilde{R}_{2,2} = \tilde{R}_{4,4} = \cdots = \tilde{R}_{n,n} = 1$, 
$\tilde{Q}_{1,1} = 4, \tilde{Q}_{2,2} = \tilde{Q}_{4,4} = \cdots = \tilde{Q}_{n,n} = \frac{1}{4}$,
$\tilde{Q}_{3,3} = \tilde{Q}_{5,5} = \cdots = \tilde{Q}_{(n-1),(n-1)} = 8$,
$\tilde{Q}_{1,3} = \tilde{Q}_{3,5} = \cdots = \tilde{Q}_{(n-3),(n-1)} = -4$,
$\tilde{Q}_{3,1} = \tilde{Q}_{5,3} = \cdots = \tilde{Q}_{(n-1),(n-3)} = -4$,
$\tilde{B}_{2,1} = \tilde{B}_{4,2} = 1, \tilde{C}_{1,1} = 1$,
where the other entries of $\tilde{J}$, $\tilde{R}$, $\tilde{Q}$, $\tilde{B}$, and $\tilde{C}$ are zeros.

We reduced the dimension $n=50$ to $r=4, 6, 8, 10, 30$, and compared the BT results with those from our proposed method.
Tables \ref{table1}, \ref{table2}, and \ref{table3} present the results for the $H^2$ error norm, $H^{\infty}$ error norm, and gradient norm, respectively.
In Appendix \ref{apeC}, we give the reduced matrices $(J_r, R_r,B_r, C_r)$ in the case where $r=4$. 
For each $r\in \{4, 6, 8, 10\}$, the $H^2$ and $H^{\infty}$ error norms given by the BT method are greater than those of our method.
In particular, for each $r\in \{4, 6, 8\}$, the $H^2$ error norms of  our method are less than $1/7$ of the corresponding error norms of the BT method.
This is because the reduced models of the BT method are far from optimal, as can be seen from Table \ref{table3}.
Moreover, we can conclude that our proposed method gives a globally near-optimal solution to Problem 5. 
Let $\sigma_1\geq \sigma_2 \geq \cdots \geq \sigma_r\geq \sigma_{r+1}\geq \cdots \geq \sigma_n > 0$ be the Hankel singular values associated with the realization $(A,B,C)$ of the transfer function $G$. Then, for any reduced transfer function $G_r$ of order $r$, we have that
\begin{align*}
||G-G_r||_{H^{\infty}} \geq \sigma_{r+1},
\end{align*}
as shown in \cite{dullerud2000course}.
Thus, from Table \ref{table2}, we can conclude that our proposed method provides a globally near-optimal solution to Problem 5.
Furthermore, when $r=30$, the BT method provides a locally optimal solution to Problem 2, because $||{\rm grad} f(J_r,R_r,B_r,C_r)||_{\eqref{Riemannian_metric}}$ is sufficiently close to zero,
where $||\cdot||_{\eqref{Riemannian_metric}}$ denotes the induced norm from Riemannian metric \eqref{Riemannian_metric}.

From the above observations, it can be seen that the BT method may provide locally optimal solutions to Problem 2 if the reduced model dimension is sufficiently large.
However, from the viewpoint of controller design, it is preferable for the dimension of the state of a plant to be as small as possible.
Thus, to reduce the original stable linear system to a small-dimensional system, our proposed method will be useful for improving the results of the BT method.

\begin{figure}[t]
\begin{center}
\includegraphics[width=95mm]{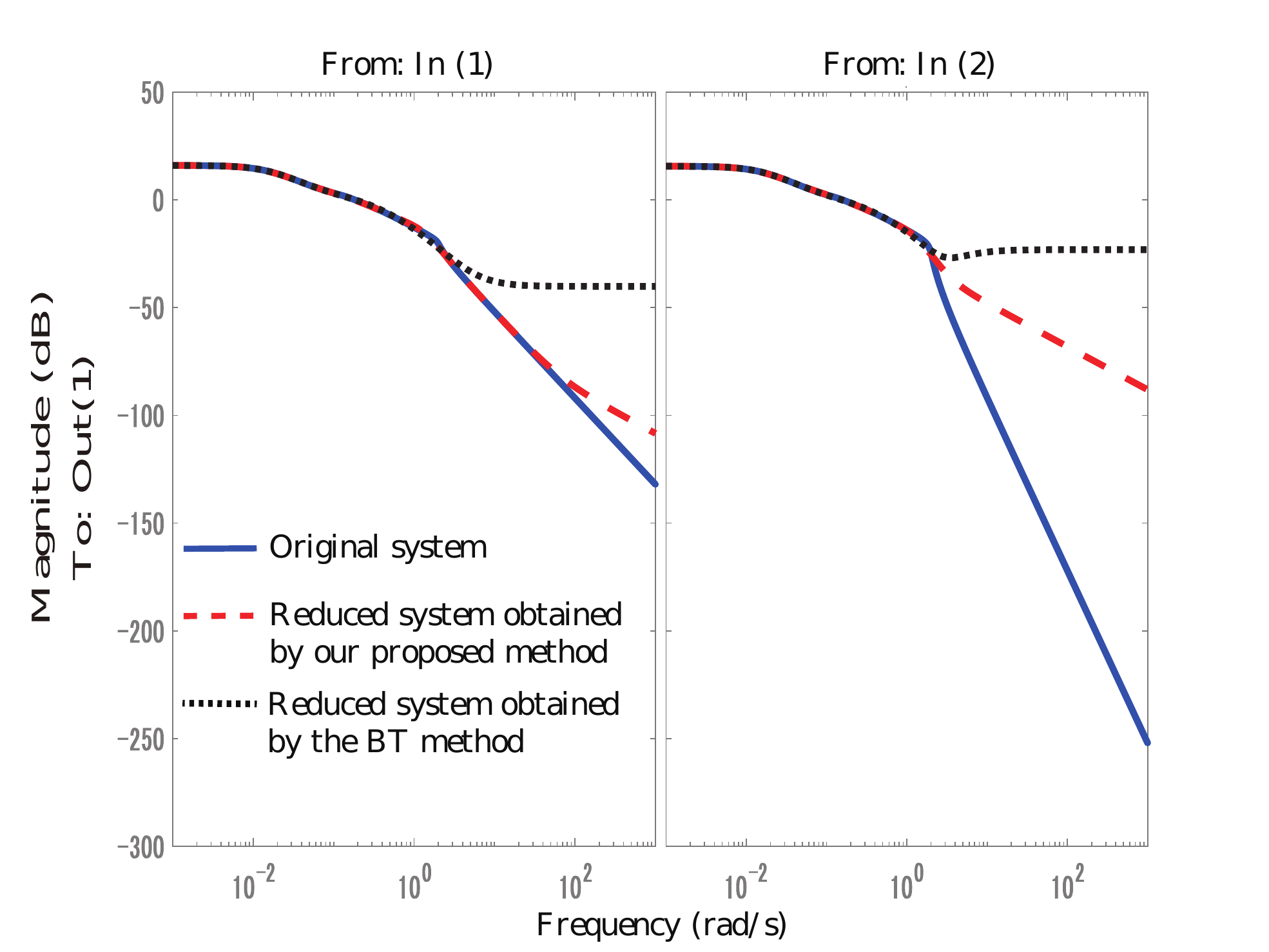}
\end{center}
\vspace{-3mm}
\caption{Bode diagram of original and reduced systems.} \label{bode}
\end{figure}

Fig.\,\ref{bode} illustrates the Bode diagram of the original system and the reduced systems
obtained by the proposed method and BT.
When the frequency is less than 1 rad/s, both reduced systems coincide with the original.
In contrast, when the frequency is greater than 1 rad/s, our reduced system is closer to the original than the system obtained by the BT method.
Thus, we can conclude that our proposed method produces better reduced systems than BT in terms of the frequency response.

\begin{figure}[t]
\begin{center}
\includegraphics[width=80mm]{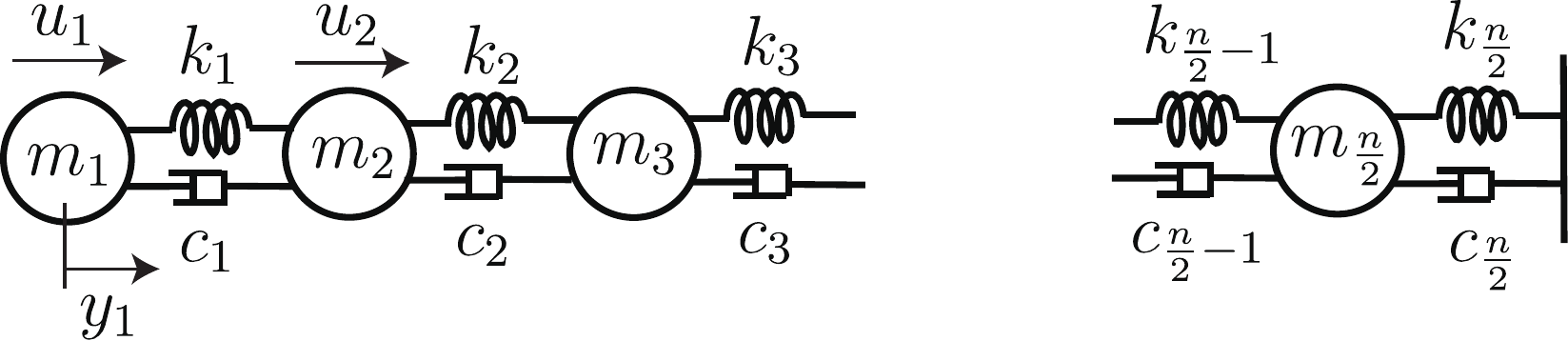}
\end{center}
\vspace{-4mm}
\caption{Mass-spring-damper system.} \label{mass}
\vspace{2mm}
\end{figure}

\begin{table*}[h]
\caption{$||G-G_r||_{H^2}$.} \label{table1}
  \begin{center}
    \begin{tabular}{|c|c|c|c|c|c|} \hline
          $r$  & 4 & 6 & 8 & 10 & 30  \\ \hline 
BT method  & 0.23248&  0.11858 & 0.05526 & 0.02416 & 0.00002  \\ 
Proposed method  & 0.03218&  0.01061 &  0.00765  &0.00552 & 0.00002  \\ \hline
    \end{tabular}
  \end{center}

\caption{$||G-G_r||_{H^{\infty}}$.} \label{table2}
  \begin{center}
    \begin{tabular}{|c|c|c|c|c|c|} \hline
          $r$  & 4 & 6 & 8 & 10 & 30 \\ \hline 
  BT method  & 0.11669& 0.05198 & 0.01646 & 0.00989 & 0.00008   \\ 
	Proposed method   & 0.04891& 0.03182 & 0.01171 & 0.00908 & 0.00008  \\
    $\sigma_{r+1}$     & 0.02834   &  0.01198 & 0.00508  &  0.00262 & 0.00004    \\ \hline
    \end{tabular}
  \end{center}

\caption{$||{\rm grad}\, f(J_r,R_r,B_r,C_r)||_{\eqref{Riemannian_metric}}$.} \label{table3}
  \begin{center}
    \begin{tabular}{|c|c|c|c|c|c|} \hline
          $r$  & 4 & 6 & 8 & 10 &30  \\ \hline 
  BT method & $4.3\times 10^{-1}$ & $1.5\times 10^{-1}$ & $5.3\times 10^{-2}$ & $1.8\times 10^{-2}$ & $3.1\times 10^{-6}$  \\ 
     Proposed method & $8.2\times 10^{-5}$ & $9.8\times 10^{-5}$ & $7.4\times 10^{-5}$ & $7.4\times 10^{-5}$ & $3.1\times 10^{-6}$     \\ \hline
    \end{tabular}
  \end{center}
\end{table*}

\subsection{Building system}

We also consider the building model of the Los Angeles University Hospital reported in \cite{chahlaoui2002collection}.
This model can be described by \eqref{1}, and has $n=48$ and $m=p=1$.
For $r=3$, our proposed method produces a reduced model with
\begin{align*}
||G-G_r||_{H_2} = 0.0030\quad {\rm and}\quad ||G-G_r||_{H^{\infty}} = 0.0039,
\end{align*}
although the BT method gives
\begin{align*}
||G-G_r||_{H_2} = 0.0416\quad {\rm and}\quad ||G-G_r||_{H^{\infty}} = 0.0079.
\end{align*}
Here, the values of $||{\rm grad}\, f(J_r,R_r,B_r,C_r)||_{\eqref{Riemannian_metric}}$ attained by BT and our method are $4.1\times 10^{-2}$ and $9.8\times 10^{-6}$, respectively.
Note that 
\begin{align*}
||G-G||_{H^{\infty}} \geq \sigma_4 =  0.0019.
\end{align*}
Thus, we conclude that our proposed method provides a globally near-optimal solution to Problem 5.
Furthermore, from the above results, our proposed method produces a better reduced system than the BT method.

\begin{figure}[t]
\begin{center}
\includegraphics[width=95mm]{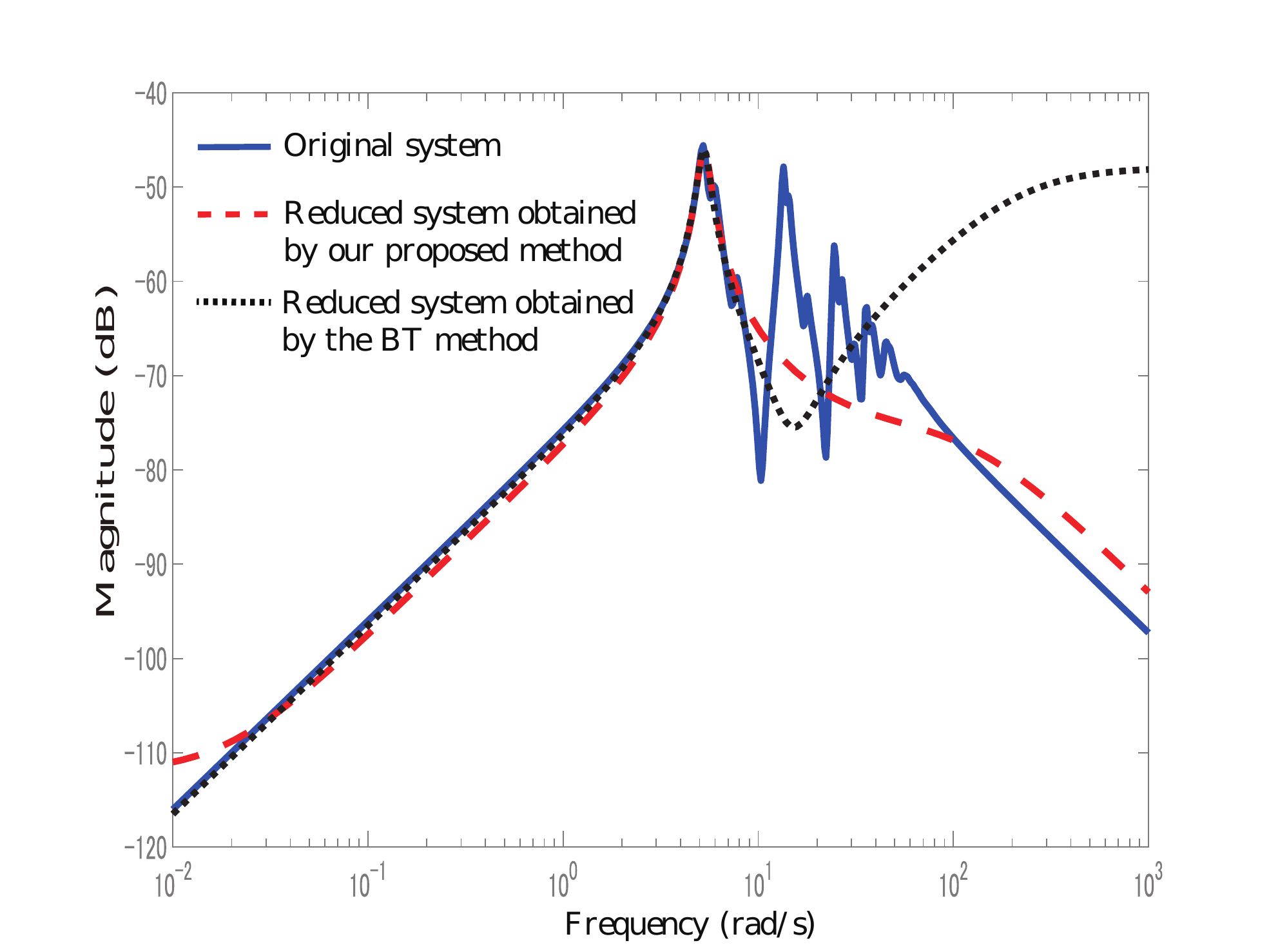}
\end{center}
\vspace{-3mm}
\caption{Bode diagram of original and reduced systems.} \label{bode2}
\end{figure}

Fig.\,\ref{bode2} illustrates the Bode diagram of the original and reduced systems.
The results given by the proposed method are similar to those of the BT methods in the low-frequency region.
However, our proposed method gives considerably better results than BT at higher frequencies.

\section{Conclusion} \label{sec6}

We have proposed a Riemannian optimal model reduction method for stable linear systems.
The model reduction problem was formulated as a minimization problem of the $H^2$ error norm between the transfer functions of the original and reduced systems
on the product manifold of the set of skew-symmetric matrices, the manifold of  the symmetric positive-definite matrices, and two Euclidean spaces.
The stability of the reduced systems constructed using the optimal solutions to our problem is preserved.
Moreover, we proposed that the initial point in our algorithm should be the output of the BT method,
because BT produces satisfactory reduced models and is easily implemented in MATLAB.
Numerical experiments demonstrated that, in the sense of the $H^2$ norm, our method achieves outstanding performance compared with the BT method when the reduced model dimension is small.
Furthermore, we illustrated that our method provides globally near-optimal solutions to the minimization problem of the $H^{\infty}$ error norm.
Moreover, Bode diagrams showed that our method is better than the BT method.

\section*{Acknowledgment}

This work was supported by JSPS KAKENHI Grant Number JP18K13773.

 \appendix


\subsection{Proof of \eqref{3_2}} \label{ape0}

For convenience, we prove \eqref{3_2}, although a similar discussion can be found in \cite{sato2017}.

Because systems \eqref{1} and \eqref{remark_red} are both asymptotically stable, they are $L^2$-stable.
That is, $u\in L^2({\bf R}^m)$ implies that $y, \hat{y}_r\in L^2({\bf R}^p)$, and thus,
there exist Fourier transformations $U$, $Y$, and $\hat{Y}_r$ of $u$, $y$, and $\hat{y}_r$, respectively.
Hence, we have that
\begin{align*}
& ||y-\hat{y}_r||_{L^{\infty}}\\
=& \sup_{t\geq 0} ||y(t)-\hat{y}_r(t)|| \\
=& \sup_{t\geq 0} ||\frac{1}{2\pi} \int_{-\infty}^{\infty} (Y(i\omega)-\hat{Y}_r(i\omega)) e^{i\omega t} d\omega || \\
\leq & \frac{1}{2\pi} \int_{-\infty}^{\infty} ||Y(i\omega)-\hat{Y}_r(i\omega)|| d\omega \\
\leq & \frac{1}{2\pi} \int_{-\infty}^{\infty} ||G(i\omega)-\hat{G}_r(i\omega)|| \cdot ||U(i\omega)|| d\omega \\
\leq & \sqrt{\frac{1}{2\pi} \int_{-\infty}^{\infty} ||G(i\omega)-\hat{G}_r(i\omega)||^2d\omega}\sqrt{\frac{1}{2\pi} \int_{-\infty}^{\infty} ||U(i\omega)||^2d\omega}\\
\leq &||G-\hat{G}_r||_{H^2}\cdot ||u||_{L^2},
\end{align*}
where the second equality follows from the inverse Fourier transformations of $Y$ and $\hat{Y}_r$, the fifth inequality is from the Cauchy--Schwarz inequality, and the
final inequality follows from $||G(i\omega)-\hat{G}_r(i\omega)||\leq ||G(i\omega)-\hat{G}_r(i\omega)||_F$ and Parseval's theorem.
This completes the proof.

\subsection{Geometry of the manifold ${\rm Sym}_+(r)$} \label{apeB}
We review the geometry of ${\rm Sym}_{+}(r)$ to develop an optimization algorithm for solving Problem 1.
For a detailed explanation, see \cite{sato2018TAC}.

For $\xi_1$, $\xi_2\in T_S {\rm Sym}_{+}(r)$, we define the Riemannian metric as
\begin{align}
\langle \xi_1, \xi_2\rangle_S := {\rm tr} ( S^{-1} \xi_1 S^{-1} \xi_2  ). \label{metric}
\end{align}
Let $g: {\rm Sym}_+(r) \rightarrow {\bf R}$ be a smooth function and $\bar{g}$ be the extension of $g$ to the Euclidean space ${\bf R}^{r\times r}$.
The Riemannian gradient ${\rm grad}\, f(S)$ with respect to the Riemannian metric \eqref{metric} is given by
\begin{align}
{\rm grad}\, g(S) = S {\rm sym}(\nabla \bar{g}(S)) S, \label{gradient}
\end{align}
where $\nabla \bar{g}(S)$ denotes the Euclidean gradient of $\bar{g}$ at $S\in {\rm Sym}_+(r)$.
The Riemannian Hessian ${\rm Hess}\,g(S): T_S {\rm Sym}_+(r) \rightarrow T_S {\rm Sym}_+(r)$ of the function $g$ at $S\in {\rm Sym}_+(r)$ is given by
\begin{align}
{\rm Hess}\,g(S)[\xi] =& S {\rm sym}( {\rm D} \nabla \bar{g}(S) [\xi] )S \nonumber\\
& +  {\rm sym}(\xi {\rm sym} (\nabla \bar{g}(S)) S ). \label{Hess}
\end{align}
The exponential map on ${\rm Sym}_+(r)$ is given by
\begin{align}
{\rm Exp}_S (\xi) = S^{\frac{1}{2}} \exp (S^{-\frac{1}{2}} \xi S^{-\frac{1}{2}} ) S^{\frac{1}{2}}, \label{8}
\end{align}
where $\exp$ is the matrix exponential function.

\subsection{Reduced matrices $(J_r,R_r,B_r,C_r)$ in the case where $r=4$ in Section \ref{sec4}} \label{apeC}

We present the reduced matrices $(J_r,R_r,B_r,C_r)$ produced by our proposed method in the case where $r=4$ in Section \ref{sec4}. 
Let
$J_r=\begin{pmatrix}
(J_r)^1 & (J_r)^2
\end{pmatrix}$, $R_r = \begin{pmatrix}
(R_r)^1 & (R_r)^2
\end{pmatrix}$,
$C_r = \begin{pmatrix}
(C_r)^1 & (C_r)^2
\end{pmatrix}$.
Then, we obtain
\begin{align*}
(J_r)^1 &= \begin{pmatrix}
0.000000000000000 &  -0.049530743507566 \\
   0.049530743507566 &    0.000000000000000  \\
  -0.018625039127746 &  0.626524211054092  \\
   0.007106890495913  & -1.083765311671058 
\end{pmatrix}, \\
(J_r)^2 &=\begin{pmatrix}
  0.018625039127746 & -0.007106890495913\\
    -0.626524211054092 &  1.083765311671058 \\
    0.000000000000000 &  0.066881602488369 \\
  -0.066881602488369 & 0.000000000000000
\end{pmatrix}, \\
(R_r)^1 &= \begin{pmatrix}
   0.020979798103068 &  0.008729495305520  \\
   0.008729495305520  & 0.296162218193050  \\
  -0.026753473825891  & 0.016509857981159  \\
  -0.003019900398660 & -0.169695898367632 
\end{pmatrix}, \\
(R_r)^2 &= \begin{pmatrix}
  -0.026753473825891 & -0.003019900398660 \\
  0.016509857981159 & -0.169695898367632\\
  0.277287705425208  &-0.447429037737505\\
  -0.447429037737505  & 1.303620534440710
\end{pmatrix}, \\
B_r & = \begin{pmatrix}
   1.087281955207546 &  1.075128712585373 \\
   0.019632883027025 & -0.081897882654859 \\
  -0.060704161404099 & -0.031902870273656 \\
   0.013609328117831 & -0.011572768539278
\end{pmatrix}, \\
(C_r)^1 & = \begin{pmatrix}
0.079020553332377 &  0.648595865888539 
\end{pmatrix}, \\
(C_r)^2 & = \begin{pmatrix}
0.877453660076422 & -3.055799879863735
\end{pmatrix}.
\end{align*}
\end{document}